\documentclass{article}

\usepackage{makeidx}
\usepackage{latexsym}
\usepackage{amsmath}
\usepackage{amssymb}
\usepackage{amscd}
\usepackage{theorem}
\usepackage{multicol}
\usepackage{color}
\usepackage[all,cmtip]{xy}
\usepackage{graphicx}
\usepackage{wrapfig}

\newtheorem{theorem}{Theorem}[section]

\newtheorem{teorema}[theorem]{Theorem}

\newtheorem{proposicion}[theorem]{Proposition}

\newtheorem{lema}[theorem]{Lemma}

\newtheorem{corolario}[theorem]{Corollary}

{\theorembodyfont{\rmfamily }
\newtheorem{remark}[theorem]{Remark}}

{\theorembodyfont{\rmfamily }
}

{\theorembodyfont{\rmfamily }
\newtheorem{nada}[theorem]{}}

{\theorembodyfont{\rmfamily }
}

\usepackage[applemac]{inputenc}

\begin{document}

\title{THE FUNDAMENTAL THEOREM FOR LOCALLY PROJECTIVE GEOMETRIES}
\author{Juan B. Sancho of Salas\smallskip\\{\small Dpto. of Matem\'{a}ticas, Univ. of Extremadura, E-06071 Badajoz (SPAIN)}\\
 {\small {\it E-mail:} jsancho@unex.es}}

\date{}

\maketitle

\smallskip

\begin{abstract}
We deal with generalizations of the Fundamental Theorem of Projective Geometry to other related geometries (of dimension $\geq 3$) and non bijective maps. We consider locally projective geometries and locally affino-projective geometries (which include the classical M\"obius, Laguerre and Minkowski geometries).\end{abstract}

{\small \noindent{\it Key words:} Fundamental Theorem, locally projective geometries, affino-projective geometries, M\"{o}bius geometry\medskip

}

\section*{Introduction}

Let $V,V'$ be vector spaces over certain division rings $K,K'$ respectively, and let $\mathbb{P}(V),\mathbb{P}(V')$ be the associated projective spaces. When these spaces have dimensions $\geq 2$, the classical Fundamental Theorem \cite{Baer} of the  Projective Geometry states that collineations $\mathbb{P}(V)\to\mathbb{P}(V')$ (bijective maps preserving lines) are just maps induced by semilinear isomorphisms $V\to V'$.

The classical theorem was generalized to non necessarily bijective maps by Faure--Fr\"olicher \cite{FF1} and Havlicek \cite{Havlicek}. Our aim is to extend this generalization to some geometries, of dimension $\geq 3$, closely related to projective spaces.

By definition, a locally projective geometry $X$ fulfills that the quotient $X/\!\!/x$ is a projective space for any point $x\in X$. Under general conditions, Wyler \cite{Wyler} in dimension $\geq 4$ and Kahn \cite{Kahn} in dimension $3$, proved that locally projective geometries admit a canonical embedding in a projective space. Under such conditions, we obtain the following generalization of the Fundamental Theorem,\medskip

\noindent\textbf{Theorem.} {\it Let $X,X'$ be locally projective geometries of dimensions $\geq 3$ and let $X\hookrightarrow\mathbb{P}(V)$, $X'\hookrightarrow\mathbb{P}(V')$ be the respective minimal projective embeddings (which exist under general conditions).

Any morphism of geometries $X\to X'$, whose image is not contained in a plane, is induced by a semilinear map $V\to V'$, unique up to a scalar factor.}\medskip

Let us resume the content of this paper.\smallskip 

In section 1 we present the category where we work. We consider closure spaces satisfying MacLane--Steinitz exchange axiom and a finitary axiom. For the sake of brevity we name them simply {\it geometries} (they are named dimensional linear spaces by Delandtsheer \cite{Buekenhout}). {\it Morphisms of geometries} are maps compatible with the closure operator. In this context, we reformulate the Fundamental Theorem for partial maps $\mathbb{P}(V)\dashrightarrow\mathbb{P}(V')$ obtained by Faurer--Fr\"olicher and Havlicek.

In section 2, using this theorem, we prove the following main technical result:\medskip

\noindent\textbf{Theorem.} {\it Let $\mathbb{P}=\mathbb{P}(V)$ and $\mathbb{P}'=\mathbb{P}(V')$ be projective spaces of dimensions $\geq 3$ over division rings $K$ and $K'$, respectively. Let $X\subseteq\mathbb{P}$ be a subgeometry with enough points and let $\varphi\colon X\to\mathbb{P}'$ be a morphism such that:\smallskip

$1$. The image of $\varphi$ is not contained in a plane of $\mathbb{P}'$.\smallskip

$2$. There exist $x_1,x_2\in X$ such that $X/\!\!/x_i=\mathbb{P}/\!\!/x_i$ ($i=1,2$) and $\varphi(x_1)\neq\varphi(x_2)$.\smallskip

Then the map $\varphi\colon X\to\mathbb{P}'$ is induced by a semilinear map $\Phi\colon V\to V'$, unique up to a scalar factor.}\medskip

As a direct consequence, in section 3 we obtain the Fundamental Theorem for locally projective geometries stated above. In this section we also recollect several elementary facts of those geometries, including an axiomatization in terms of incidence of points, lines and planes.

Finally, in section 4 we extend the Fundamental Theorem to more general geometries, named locally affino-projective geometries. A subgeometry $X$ of a projective space $\mathbb{P}$ is said to be locally affino-projective when for any point $x\in X$ there exists a hyperplane $x\in H_x\subset\mathbb{P}$ such that $\mathbb{P}/\!\!/x=(X/\!\!/x)\cup(H_x/\!\!/x)$; this condition means that any line $x\in L\subset\mathbb{P}$ fulfills $\,|L\cap X|\geq 2\,$ or $\,L\subseteq H_x$. In the case of the real projective space of dimension $3$, these geometries include the so named Benz planes: The M\"obius geometry of an ellipsoid, the Minkowski geometry of a ruled quadric and the Laguerre geometry of a cone. For these geometries we obtain the following theorem.\medskip

\noindent\textbf{Theorem.} {\it Let $\mathbb{P}=\mathbb{P}(V)$ and $\,\mathbb{P}\,'=\mathbb{P}(V')$ be projective spaces of dimensions $\geq 3$ over division rings $K$ and $K'$, respectively, with $|K|\geq 4$ or $|K|=3=\mathrm{char}\, K'$. Let $X\subseteq\mathbb{P}$ be a locally affino-projective subgeometry with enough points.

Any morphism $\varphi\colon X\to\mathbb{P}\,'$, whose image is not contained in a plane, is induced by a semilinear map $\Phi\colon V\to V'$, unique up to a scalar factor.}

\newpage

\section{Preliminares}

In this preliminary section we fix the category where we shall work.

\subsection*{Geometries}

The spaces that we shall consider have different names in the literature. For the sake of brevity we shall use the name {\it geometry}, also used in \cite{FF} and \cite{Gratzer}.\medskip

\noindent\textbf{Definition.} A \textbf{geometry} is a set $X$ (whose elements are named \textbf{points}) with a family $\mathfrak{F}$ of subsets of $X$ (named  \textbf{subspaces} or \textbf{flats}), satisfying the following four axioms:\smallskip

\noindent\textbf{G1.} $\emptyset\in\mathfrak{F}$, $X\in\mathfrak{F}$ and $\{x\}\in\mathfrak{F}$ for any $x\in X$.\smallskip

The subspace $\{ x\}$ will be simply denoted by $x$.\smallskip

\noindent\textbf{G2.} Any intersection $\bigcap S_i$ of subspaces is a subspace.\smallskip

This axiom implies that the set $\mathfrak{F}$ of subspaces of $X$, with the inclusion, is a complete lattice: For any set $\{ S_i\}$ of subspaces, the infimum is
$\bigwedge S_i=\bigcap S_i$ and the supremum is $\bigvee S_i=\bigcap T_k$, where $T_k$ runs over the subspaces of $X$ containing $\bigcup S_i$.\medskip

\noindent\textbf{G3 (exchange axiom).} Given a subspace $S$ and a point $x\notin S$, there is no subspace $S'$ such that $\, S\subset S'\subset S\vee x\,$ (strict inclusions),\smallskip

That is to say, if $x,y\notin S$ then $\quad y\in S\vee x\,\Leftrightarrow\, x\in S\vee y$.\smallskip

The \textbf{closure} of a subset $A\subseteq X$ is the least subspace $\overline A$ containing $A$.\smallskip

\noindent\textbf{G4 (finitary axiom).} For any subset $A\subseteq X$ and any point $x\in \overline A$, we have $x\in \overline{\{a_1,\dots,a_r\}}$ for some finite subset  $\{a_1,\dots,a_r\}\subseteq A$.\medskip

% A geometry is a particular case of a \textbf{closure space:} A set $X$ endowed with a map $\mathcal{P}(X)\to\mathcal{P}(X)$, $A\mapsto\overline A$, satisfying the axioms, $$(i)\quad A\subseteq\overline A\qquad,\qquad (ii)\quad A\subseteq\overline B\quad\Rightarrow\quad\overline A\subseteq\overline B$$

%\smallskip

\noindent\textbf{Definitions.} Given a geometry $X$, a subset $A\subseteq X$ is said to be \textbf{independent} when for any $a\in A$ we have $a\notin\overline{A\setminus a}$.

Let $A$ be a subset of a subspace $S$. If $\,\overline A=S$ we say that $A$ \textbf{generates} $S$. If $A$ is independent and generates $S$, we say that $A$ is a \textbf{basis} of $S$.

\begin{theorem} Let $S$ be a subspace of a geometry $X$. Then:\smallskip

-- Any set generating $S$ contains a basis of $S$.\smallskip

-- Any independent subset of $S$ is contained in a basis of $S$.\smallskip

-- Any two bases of $S$ have equal cardinal.
\end{theorem}

For a proof in a more general context, see \cite{FF} Chap. IV.\medskip

\noindent\textbf{Definition.} The \textbf{dimension} of a subspace is the common cardinal of its bases minus 1.\medskip

Points are just subspaces of dimension $0$. Subspaces of dimension $1$ are named \textbf{lines} and subspaces of dimension $2$, \textbf{planes}.   

 Let $A=\{p_0,\dots,p_r\}\subseteq X$ be an independent finite set. Then $\overline A=p_0\vee\cdots\vee p_r$ is a subspace of dimension $r$.

In particular, any two distinct points $p_0,p_1$ lie in a unique line $\,p_0\vee p_1$, and any three non collinear points $p_0,p_1,p_2$ lie in a unique plane $\,p_0\vee p_1\vee p_2$.

\subsection*{Morphisms of geometries}

We shall come across partially defined maps. The notation $\varphi\colon X\dashrightarrow X'$ is used in order to point out that $\varphi$ is a map with domain $\text{dom}\,\varphi\subseteq X$ and image set $\text{im}\,\varphi\subseteq X'$. The \textbf{exceptional set} of $\varphi$ (with respect to $X$) is $E:=X-\text{dom}\,\varphi$. The map $\varphi$ is said to be globally defined (with respect to $X$) provided that $E$ is empty. The notation $\varphi\colon X\to X'$ is only maintained for globally defined maps.

Given any subset $S\subseteq X$, we define $\varphi(S):=\{\varphi(s):s\in S\cap\text{dom}\,\varphi\}$.  Hence $\varphi(S)=\emptyset$ exactly for $S\subseteq E$.
Given any subset $S'\subseteq X'$, we define $\varphi^{-1}(S'):=\{x\in X:\varphi(x)\in S'\text{ or }\varphi(x)=\emptyset\}=\{x\in\text{dom}\,\varphi:\varphi(x)\in S'\}\cup E$.\medskip

\noindent\textbf{Definition.} A partial map $\varphi\colon X\dashrightarrow X'$, between geometries, is a \textbf{partial morphism} when it fulfills the following equivalent conditions:\smallskip

a) The preimage of any subspace $S'\subseteq X'$ is a subspace $\varphi^{-1}(S')\subseteq X$.\smallskip

b) For any subset $A\subseteq X$ we have $\,\varphi(\,\overline A)\subseteq\overline{\varphi(A)}$, that is to say,
$$x\,\in\,\overline A\quad\Rightarrow\quad\varphi(x)\,\subseteq\,\overline{\varphi(A)}\quad .$$

c) For any finite subset $A\subseteq X$ we have $\,\varphi(\,\overline A)\subseteq\overline{\varphi(A)}$.\medskip

It is easy to check the equivalence of these conditions (see \cite{Sancho}, Prop. 1.2).\medskip

If a partial morphism is globally defined then it is named \textbf{morphism of geometries}.

\begin{nada}\label{1.4a} According $(c)$, a partial map $\varphi\colon X\dashrightarrow X'$, between geometries, is a partial morphism if and only if for any finite sequence  $\,x_0,\dots,x_r\in X\,$ we have
$$x_0\,\in\, x_1\vee\dots\vee x_r\qquad\Rightarrow\qquad \varphi(x_0)\,\subseteq\,\varphi(x_1)\vee\dots\vee\varphi(x_r)\quad .$$
Recall that $\varphi(x)$ can be the empty set.
\end{nada}

The exceptional set $E=\varphi^{-1}(\emptyset)$ of a partial morphism $\varphi\colon X\dashrightarrow X'$ is a subspace of $X$.\medskip

\noindent\textbf{Definition.} We say that a geometry $X$ is \textbf{generated by lines} when the subspaces of $X$ are just the subsets $S\subseteq X$ such that
$$x_1,x_2\in S\quad \Rightarrow\quad  x_1\vee x_2\,\subseteq S$$
That is to say, in a geometry generated by lines, subspaces are just subsets containing the line joining any two different points of such subset.\smallskip

We say that a geometry $X$ is \textbf{generated by lines and planes} when the subspaces of $X$ are just the subsets $S\subseteq X$ such that
$$x_1,x_2,x_3\in S\quad\Rightarrow\quad x_1\vee x_2\vee x_3\,\subseteq S$$
That is to say, a subset $S$ is a subspace if it contains the line joining any two different points of $S$ and the plane defined by any three non collinear points of $S$.\medskip

Projective spaces are geometries generated by lines. Affine spaces over a division ring $K\neq\mathbb{Z}_2$ also are generated by lines, and affine spaces over $\mathbb{Z}_2$ are generated by lines and planes.

\begin{nada}\label{1.6} Let $\varphi\colon X\dashrightarrow X'$ be a partial map between geometries, where $X$ is generated by lines. The map $\varphi$ is a partial morphism if and only if
\begin{equation}\label{eq}
x_0\,\in\, x_1\vee x_2\qquad\Rightarrow\qquad\varphi(x_0)\,\subseteq\,\varphi(x_1)\vee\varphi(x_2)
\end{equation}
for any $x_0,x_1,x_2\in X$. 

From (\ref{eq}) it easily follows (when $X$ is generated by lines) that the preimage of any subspace also is a subspace, hence $\varphi$ is a morphism. The converse holds by  \ref{1.4a}.\smallskip

%Condition (\ref{eq}) is equivalent to the following three conditions: 

%--  $\varphi$ transforms any three collinear points $x_0,x_1,x_2\in\text{dom}\,\varphi$ into collinear points; 

%--  for any line $L\nsubseteq E\,$ the restriction $\,\varphi_{|L\!-\!E}\,$  is injective or constant;

 %-- if $e\in E$ then $\,\varphi(e\vee x)=\varphi(x)\,$ for all $x\in X$.\medskip
 
 Partial maps satisfying  (\ref{eq}) are named {\it weak linear mappings} in \cite{Havlicek2}.
\end{nada}

\begin{nada}\label{1.7} Let $\varphi\colon X\dashrightarrow X'$ be a partial map between geometries, where $X$ is generated by lines and planes. The map $\varphi$ is a partial morphism if and only if
$$x_0\,\in\, x_1\vee x_2\vee x_3\qquad\Rightarrow\qquad\varphi(x_0)\,\subseteq\,\varphi(x_1)\vee\varphi(x_2)\vee\varphi(x_3)$$
for any $x_0,x_1,x_2,x_3\in X$.\end{nada}

\noindent\textbf{Definition.} A morphism of geometries $\varphi\colon X\to X'$ is an \textbf{isomorphism} if it is bijective and the inverse map $\varphi^{-1}$ also is a  morphism.\medskip

A bijective morphism may be not an isomorphism. For example, given a geometry $X$ of dimension $n\geq 2$ and a natural number $m<n$, let $X'$ be the geometry (named {\it  truncation} of $X$) with the same underlying set $X$ but such that the proper subspaces are just the subspaces of $X$ of dimension $<m$ (so that $\text{dim}\, X'=m$). The identity $X\to X'$ is a bijective morphism but not an isomorphism. There are more interesting examples; Ceccherini \cite{C} gives a bijective morphism between a 4-dimensional projective space and a non-arguesian projective plane.\medskip

An isomorphism $\varphi\colon X\to X'$, between geometries generated by lines, is also said to be a \textbf{collineation}. That is to say, a collineation (between such geometries) is just a bijection transforming lines onto lines.

\begin{nada}\label{1.5a}The following two statements are easy to prove (see \cite{Sancho} Prop. 1.5 and, in the context of {\it linear spaces}, \cite{Kreuzer}):\smallskip

{\it a)  If $\varphi\colon X\to X'$ is a surjective morphism of geometries then $\,\mathrm{dim}\,X\,\geq\,\mathrm{dim}\, X'$.\smallskip

b) If $\varphi\colon X\to X'$ is a surjective morphism of geometries and $\,\mathrm{dim}\,X=\mathrm{dim}\, X'<\infty\,$ then $\varphi$ is an isomorphism.}
\end{nada}

\noindent\textbf{Definition.} Given a geometry $X$, any subset $A\subseteq X$ is a geometry where subspaces are just subsets $S\cap A$ with $S$ a  subspace of $X$. Then we say that $A$ is a \textbf{subgeometry} of $X$.

The natural inclusion $A\hookrightarrow X$ is a morphism.\goodbreak\medskip

\noindent\textbf{Definition.} A morphism $\varphi\colon X\to X'$ is an \textbf{embedding} when $\varphi\colon X\to\varphi(X)$ is an isomorphism, where $\varphi(X)$ is considered as a subgeometry of $X'$.

\begin{proposicion}\label{1.6a} Let $\varphi\colon X\to X'$ be a morphism of geometries. Let $\{X_i\}_{i\in I}$ be a directed family of subspaces of $X$ such that $X=\bigcup X_i$. If $\varphi:X_i\to X'$ is an embedding for any $X_i$, then $\varphi\colon X\to X'$ is an embedding.
\end{proposicion}

The condition that the family $\{X_i\}_{i\in I}$ is directed means that for any pair of indices $i,j\in I$ there is $k\in I$ such that $X_i,X_j\subseteq X_k$.\medskip

\noindent{\it Proof.} Since $\varphi\colon X_i\to\varphi(X_i)$ is bijective for any $X_i$ and the family is directed, it follows that the morphism $\,\varphi\colon X=\bigcup X_i\longrightarrow\bigcup\varphi(X_i)=\varphi(X)\,$ is bijective. We have to see that $\varphi^{–1}\colon\varphi(X)\to X$ is a morphism. Given a finite subset $\{\varphi(x_0),\dots,\varphi(x_r)\}\subseteq\varphi(X)$ such that $\varphi(x_0)\in\varphi(x_1)\vee\cdots\vee\varphi(x_r)$, we must prove that $x_0\in x_1\vee\cdots\vee x_r$. We have $x_1,\dots, x_r\in X_i$ for some $X_i$, hence $\varphi(x_1),\dots,\varphi(x_r)\in \varphi(X_i)$; since $\varphi^{-1}\colon \varphi(X_i)\to X_i$ is a morphism we conclude that $x_0\in x_1\vee\cdots\vee x_r$.

\hfill$\square$

\subsection*{Quotients}

Partial morphisms are reduced to global morphisms due to the following notion.\medskip

\noindent\textbf{Definition.} Let $E$ be a subspace of a geometry $X$. The \textbf{quotient space} is the quotient set
$$X/\!\!/E:=\,(X\!-\! E)/\equiv$$
where $\,\equiv\,$ is the following equivalence relation on $X-E$ :
$$x_1\equiv\,x_2\quad\Leftrightarrow\quad x_1\vee E\,=\,x_2\vee E\quad .$$

Let $\,\pi\colon X\dashrightarrow X/\!\!/E$, $x\mapsto [x]$, be the quotient (partial) map. We name \textbf{subspace} of $X/\!\!/E$ any subset $\,S/\!\!/E=\pi(S)\,$ where $S$ is a subspace of $X$ containing $E$. The quotient $X/\!\!/E$, endowed with the lattice of its subspaces, is a geometry. The quotient map $\,\pi\colon X\dashrightarrow X/\!\!/E$ is a partial morphism with exceptional subspace $E$, since $\pi^{-1}(S/\!\!/E)=S$.

The quotient space has the expected universal property (see \cite{Sancho} Prop. 1.9):

\begin{nada}\label{5.9} 
 Any partial morphism $\varphi\colon X\dashrightarrow X'$, with exceptional subspace $E$, uniquely factors $\varphi=\widetilde\varphi\circ\pi$ through a morphism $\,\widetilde\varphi\colon X/\!\!/E\longrightarrow X'$.
 $$\xymatrixcolsep{1pc}\xymatrixrowsep{2pc}\xymatrix{
X  \ar@{-->}[rd]_{\varphi} \ar@{-->}[rr]^{\pi} &  & X/\!\!/E \ar[ld]^{\widetilde\varphi} \\
 & X'
}$$
\end{nada}\goodbreak

\begin{proposicion}\label{2.1b} Let $\varphi\colon X\to X'$ be a morphism of geometries.\smallskip

a) Given subspaces $E'\subseteq X'$, $E=\varphi^{-1}(E')$, the map
$$\begin{CD}
\varphi_E:X/\!\!/E @>>> X'/\!\!/E'\qquad, \qquad [x]\longmapsto[\varphi(x)]
\end{CD}$$
is a morphism of geometries.\smallskip

b) Given $x_0\in X$, $x_0'=\varphi(x_0)$, $F=\varphi^{-1}(x_0')$, the map
$$\begin{CD}
\varphi_{x_0}:X/\!\!/x_0 &\quad  --\!\to\quad & X'/\!\!/x_0'\qquad ,\qquad [x]\longmapsto [\varphi(x)]
\end{CD}$$
is a partial morphism of geometries with exceptional subspace $F/\!\!/x_0$.
\end{proposicion}
{\it Proof.} a). By \ref{5.9}, the composition partial morphism $X\overset{\varphi}{\longrightarrow}X'\overset{\pi}{\dashrightarrow}X'/\!\!/E'$, with exceptional subspace $E$, induces a morphism $\varphi_E=\widetilde{\pi\varphi}\colon X/\!\!/E\to X'/\!\!/E'$.\smallskip

b). $\varphi_{x_0}$ is the composition partial morphism $$X/\!\!/x_0\,\dashrightarrow\,(X/\!\!/x_0)/\!\!/(F/\!\!/x_0)=X/\!\!/F\,\overset{\varphi_F}{\longrightarrow}\, X'/\!\!/x_0'\quad .$$
\hfill$\square$

\begin{corolario}\label{2.3} Let $i\colon X\hookrightarrow X'$ a subgeometry. Let $E\subseteq X$ be a subspace and let $E'$ be the closure of $E$ in $X'$. Then $\,i_E\colon X/\!\!/E\to X'/\!\!/E'\,$ is a subgeometry.

In particular, for any $x\in X$ one has that $X/\!\!/x$ is a subgeometry of $X'/\!\!/x$.
\end{corolario}
{\it Proof.}  If $S'$ denotes the closure in $X'$ of a subspace $S$ of $X$, 
it is easy to check that $\,S=i^{-1}S'=S'\cap X$. The morphism $i_E\colon X/\!\!/E\hookrightarrow X'/\!\!/E'$, $[x]\mapsto[x]$, is an embedding because any subspace $\,S/\!\!/E\,$ of $\,X/\!\!/E\,$ is the inverse image by $i_E$ of the subspace $\,S'/\!\!/E'\,$ of $\,X'/\!\!/E'$,
$$i_E^{-1}(S'/\!\!/E')\,=\, (i^{-1}S')/\!\!/E\,=\, S/\!\!/E\quad .$$
\hfill$\square$

\subsection*{Projective spaces}

\textbf{Definition.} A \textbf{projective space} (eventually reducible) is a set $\mathbb{P}$ (whose elements are named \textbf{points}) with a family $\mathfrak{L}$ of subsets (named \textbf{lines}), satisfying the following three axioms:\smallskip

\textbf{P1.} For any two different points passes a unique line.\smallskip

\textbf{P2.} Any line has two different points.
\smallskip

\textbf{P3 (Veblen--Young axiom).} If a line meets two sides of a triangle, not in the common vertex, then it also meets the third side.\medskip

A \textbf{subspace} of $\mathbb{P}$ is a subset $S$ such that it contains the line joining any two different points of $S$.\medskip

It is clear that any subspace of a projective space also is a projective space.\medskip

The following two statements are well known (see for example \cite{FF}):

\begin{nada} {\it Any projective space, endowed with the lattice of subspaces, is a geometry (obviously generated by lines).}
\end{nada}

\begin{nada}\label{1.10} {\it A geometry $X$ is a projective space if and only if for any two subspaces $S_1,S_2\subseteq X$ of finite dimension the} \textbf{dimension formula} {\it holds},
$$\text{dim}\,S_1+\text{dim}\, S_2\,=\, \text{dim}\,(S_1\vee S_2)+\text{dim}\, (S_1\cap S_2)\quad .$$
\end{nada}

\noindent\textbf{Definition.} An  \textbf{irreducible projective space} is a projective space where any line contains at least three different points.

\begin{nada}\label{2.2a} If a projective space is not irreducible, it is named \textbf{reducible}. The following consideration explains the name.

A projective space $\mathbb{P}$ is a union (in fact, a coproduct) of irreducible projective spaces as follows. Let us consider in $\mathbb{P}$ the following equivalence relation: Two  points $p_1,p_2\in \mathbb{P}$ are equivalent if either $p_1=p_2$, or $p_1\neq p_2$ and the line $p_1\vee p_2$ at least has three points. It is easy to check that it is an equivalence relation. Then $\mathbb{P}$ is the union of the equivalence classes, each one being an irreducible projective space.
\end{nada}

\noindent\textbf{Definition.} Let $V$ be a (left)  vector space over a division ring. We name \textbf{projectivization} of $V$, or \textbf{projective space associated} to $V$, the set
$$\mathbb{P}(V):=\,\{\text{1-dimensional vector subspaces }\langle v\rangle\text{ of }V\}\quad .$$

The set $\mathbb{P}(V)$ is an irreducible projective space when the lines are considered to be the projectivization of the 2-dimensional vector subspaces of $V$.

Moreover, the subspaces of $\mathbb{P}(V)$ are just the projectivization $\mathbb{P}(W)$ of the vector subspaces $W\subseteq V$. One has $\text{dim}\,\mathbb{P}(W)=\text{dim}\, W-1$. \medskip

It is a classical result \cite{Baer} that any irreducible projective space of dimension $\geq 3$ is the projectivization of a vector space:

\begin{teorema}\label{2.17}  Any irreducible projective space, of dimension $\geq 3$ or of dimension $2$ and arguesian, is isomorphic (a collineation) to the projectivization $\mathbb{P}(V)$ of a vector space $V$ over some division ring $K$.
\end{teorema}

\noindent\textbf{Definitions.} Let $V$, $V'$ be vector spaces over division rings $K$, $K'$, respectively. A map $\varphi\colon V\to V'$ is said to be \textbf{semilinear} when:\smallskip

1) it is additive: $\varphi(v_1+v_2)=\varphi(v_1)+\varphi(v_2)\quad$ for all $v_1,v_2\in V$;\smallskip

2) there is a ring morphism $\sigma\colon K\to K'$ such that $\,\varphi(\lambda v)=\sigma(\lambda)\varphi(v)\, $ for all $\lambda\in K$, $v\in V$.\smallskip

Remark that we do not require that the morphism $\sigma\colon K\to K'$ be surjective. If $\varphi$ is not the null map, then the morphism $\sigma$ is unique and it is said to be the \textbf{ring morphism associated} to $\varphi$.\smallskip

A semilinear map $\varphi\colon V\to V'$ is said to be a \textbf{semilinear isomorphism} when $\varphi$ and $\sigma$ are bijective. In such case the inverse map  $\varphi^{-1}\colon V'\to V$ also is semilinear.\medskip

\noindent\textbf{Definitions.} A semilinear map $\Phi\colon V\to V'$, $v\mapsto v'$, induces a partial map
$$\begin{CD}
\phi :\mathbb{P}(V) &\quad ---\to\quad & \mathbb{P}(V')\qquad ,\qquad \langle v\rangle\,\mapsto\,\langle v'\rangle
\end{CD}$$
where the exceptional subspace is $E=\mathbb{P}(\text{ker}\,\Phi)$.
Then $\phi$ is named \textbf{partial projective morphism} associated to $\Phi$. 
It is easy to check that $\phi$ is a partial morphism of geometries (\cite{Sancho} Prop. 2.9).

Remark that $\phi$ is globally defined ($E=\emptyset$) if and only if $\Phi$ is injective ($\text{ker}\,\Phi=0$). In such case, we say that $\phi\colon \mathbb{P}(V)\to\mathbb{P}(V')$ is a \textbf{projective morphism}. In general a projective morphism is not injective nor surjective.

A \textbf{projective isomorphism} $\phi\colon \mathbb{P}(V)\to\mathbb{P}(V')$ is the projectivization of a semilinear isomorphism $\Phi\colon V\to V'$.\medskip

It is unknown if each bijective projective morphism $\mathbb{P}(V)\to\mathbb{P}(V')$, with dimensions $\geq 2$, is an isomorphism. There are examples \cite{Kreuzer2} of injective morphisms $\mathbb{P}(V)\to\mathbb{P}(V')$, preserving non-collinearity, with $\text{dim}\,\mathbb{P}(V)>\text{dim}\,\mathbb{P}(V')$.

\begin{nada}\label{1.13} Given a subspace  $\mathbb{P}(W)\subseteq\mathbb{P}(V)$, one has an isomorphism of geometries 
$$\mathbb{P}(V)/\!\!/\mathbb{P}(W)\,=\,\mathbb{P}(V/W)\qquad ,\qquad [\langle v\rangle ]\,\mapsto\,\langle[v]\rangle\quad .$$
Moreover, the quotient morphism $\,\pi\colon\mathbb{P}(V)\dashrightarrow \mathbb{P}(V)/\!\!/\mathbb{P}(W)=\mathbb{P}(V/W)\,$ is the partial projective morphism associated to the quotient linear map $\,V\to V/W$.
\end{nada}

\begin{nada}\label{1.14} If two semilinear maps $\Phi,\Phi'\colon V\to V'$ induce the same partial projective morphism $\phi=\phi'\colon\mathbb{P}(V)\dashrightarrow\mathbb{P}(V')$ and this is not constant, then it is easy to check that $\Phi$ and $\Phi'$ are proportional: $\Phi'=\lambda\Phi$ for some $\lambda\in K'\!-\!\{0\}$ (see for example, \cite{F} Lemma 2.4).
\end{nada}

\begin{nada}\label{2.5}\textbf{Fundamental Theorem of Projective Geometry (\cite{FF1}, \cite{Havlicek})}
{\it Let $\mathbb{P}=\mathbb{P}(V)$ and $\mathbb{P}\,'=\mathbb{P}(V')$  be projective spaces of dimensions $\geq 2$ over division rings $K$ and $K'$, respectively. Let $\phi\colon\mathbb{P}\dashrightarrow\mathbb{P}'$ be a partial map, such that the image is not contained in a line. Then  $\phi$ is a partial projective morphism if and only if $\phi$ is a partial morphism of geometries.
}
\end{nada}

This statement is a generalization of the classical theorem (which only considered bijective maps) due to Faure--Fr\"olicher \cite{FF1}  and Havlicek \cite{Havlicek}; see also \cite{Havlicek2} Th. 2. A simple proof may be seen in \cite{F}.

\section{Main theorem}

\begin{lema}\label{3.5} Let $\varphi\colon X\to X'$ be a morphism of geometries such that the image is not contained in a $r$-dimensional subspace of $X'$. Then for any $x_0\in X$, $x'_0=\varphi(x_0)$, the  image of $\varphi_{x_0}\colon X/\!\!/x_0\dashrightarrow X'/\!\!/x_0'$ is not contained in a $(r-1)$-dimensional subspace of $X'/\!\!/x_0'$.
\end{lema}
{\it Proof.}  Any $(r\!-\!1)$-subspace of $X'/\!\!/x'_0$ is $S'/\!\!/x'_0$ for some  $r$-subspace $S'$ of $X'$ passing through $x_0'$. If the image of $\varphi_{x_0}$ is contained in some $(r\!-\!1)$-subspace $S'/\!\!/x_0'$, that is to say, 
$$X/\!\!/x_0\,=\,\varphi_{x_0}^{-1}(S'/\!\!/x_0')\,=\,(\varphi^{-1}S')/\!\!/x_0\quad ,$$
then taking inverse image by $\pi\colon X\dashrightarrow X/\!\!/x_0$ we obtain $X=\varphi^{-1}S'$ against the hypothesis.

\hfill$\square$\medskip

Given two $K$-linear maps $f_1\colon V_1\to V_0$, $f_2\colon V_2\to V_0$, we may consider the $K$-vector space \textbf{fiber product}
$$V_1\times_{V_0}V_2:=\,\{(v_1,v_2)\in V_1\times V_2:f_1(v_1)=f_2(v_2)\}\quad .$$

\begin{lema}\label{3.1} Let $V$ be a $K$-vector space and let $\,V_1,V_2\subseteq V\,$ be vector subspaces such that $V_1\cap V_2=0$. Let us consider the quotient vector spaces $\,\overline V_i:=V/V_i\,$  ($i=1,2$) and $\,\overline V_{12}:=V/(V_1\oplus V_2)$. We have a $K$-linear isomorphism
$$\begin{CD}
 V @= \overline V_1\times_{\overline V_{12}}\overline V_2\qquad ,\qquad v\longmapsto ([v],[v])\quad .
\end{CD}$$
\end{lema}
{\it Proof.} The map is injective because the kernel is $V_1\cap V_2=0$. Let us see that it is also surjective: Let $([w_1],[w_2])\in\overline V_1\times_{\overline V_{12}}\overline V_2$; since the classes $[w_1],[w_2]$ coincide in $\overline V_{12}$, one has $w_1\equiv w_2\ \text{mod}\, V_1\oplus V_2$, that is to say, $w_1-w_2=v_1+v_2$ where $v_i\in V_i$, so that $w_1-v_1=w_2+v_2=:v$. These equalities show that $([v],[v])=([w_1],[w_2])$.

\hfill$\square$

\begin{nada}\label{3.4} With the notations of the lemma, the quotient maps $V\to\overline V_1$, $V\to\overline V_2$ correspond, via the isomorphism $\,V=\overline V_1\times_{\overline V_{12}}\overline V_2$, with the natural projections $\,\overline V_1\times_{\overline V_{12}}\overline V_2\longrightarrow\overline V_1$, $\,\overline V_1\times_{\overline V_{12}}\overline V_2\longrightarrow\overline V_2$. Therefore, there is a canonical isomorphism between the corresponding kernels
$$\begin{CD}
V_1 @= 0\times_{\overline V_{12}}\overline V_2 & \qquad ,\qquad & V_2 @= \overline V_1\times_{\overline V_{12}}0\quad .
\end{CD}$$
\end{nada}

\noindent\textbf{Definition.} Let $X$ be a geometry of dimension $\geq 2$. We say that $X$ has \textbf{enough points} when any plane $P\subseteq X$ has four points, no three of which are collinear. That is to say, the lines of $X/\!\!/x$ have at least three points for all $x\in X$.\smallskip

If the lines of $X$ have at least three points then it is easy to see that $X$ has enough points.

\begin{teorema}\label{3.6} Let $\mathbb{P}=\mathbb{P}(V)$ and $\mathbb{P}'=\mathbb{P}(V')$ be projective spaces of dimensions $\geq 3$ over division rings $K$ and $K'$, respectively. Let $X\subseteq\mathbb{P}$ be a subgeometry with enough points and let $\varphi\colon X\to\mathbb{P}'$, $x\mapsto x'$, be a morphism such that:\smallskip

$1)$ the image of $\varphi$ is not contained in a plane of $\,\mathbb{P}'$;\smallskip

$2)$ there exist $x_1,x_2\in X$ such that $X/\!\!/x_i=\mathbb{P}/\!\!/x_i$ ($i=1,2$) and $x_1'\neq x_2'$.\smallskip

Then $\varphi\colon X\to\mathbb{P}'$ is the restriction of a unique partial projective morphism $\,\phi\colon\mathbb{P}\dashrightarrow\mathbb{P}'$.
\end{teorema}
{\it Proof.} %Condition $1$ implies that $\text{dim}\,\varphi(X)\geq 3$, hence $\text{dim}\, X\geq 3$ by \ref{1.5a}.
 Let us consider the lines $\,L:=x_1\vee x_2\subset\mathbb{P}\,$ and $\,L':=x_1'\vee x_2'\subset\mathbb{P}'$. Put $\,L_0:=L\cap X$, which is a line of the geometry $X$.

We have $x_i=\mathbb{P}\langle v_i\rangle$ for some 1-dimensional vector subspaces $\langle v_i\rangle \subset V$ and, analogously, $x_i'=\mathbb{P}\langle v_i'\rangle$ with $\langle v_i'\rangle \subset V'$. Then $\,L=\mathbb{P}\langle v_1,v_2\rangle\,$ and $\,L'=\mathbb{P}\langle v_1',v_2'\rangle$.

We have the following commutative diagram of partial morphisms,
$$\xymatrixcolsep{5pc}\xymatrixrowsep{2pc}\xymatrix{
\mathbb{P}(V/\langle v_1\rangle)\overset{\ref{1.13}}{=}\mathbb{P}/\!\!/x_1=X/\!\!/x_1\quad \ar@{-->}[r]^{\varphi_{x_1}} \ar@{-->}[d] & \quad\mathbb{P}'/\!\!/x_1'=\mathbb{P}(V'/\langle v_1'\rangle)  \ar@{-->}[d]\\
\mathbb{P}(V/\langle v_1,v_2\rangle)=\mathbb{P}/\!\!/L=X/\!\!/L_0\quad \ar@{-->}[r]^{\varphi_{L_0}}  & \quad\mathbb{P}'/\!\!/L'=\mathbb{P}(V'/\langle v_1',v_2'\rangle)  \\
\mathbb{P}(V/\langle v_2\rangle)=\mathbb{P}/\!\!/x_2=X/\!\!/x_2\quad \ar@{-->}[r]^{\varphi_{x_2}} \ar@{-->}[u] & \quad\mathbb{P}'/\!\!/x_2'=\mathbb{P}(V'/\langle v_2'\rangle)  \ar@{-->}[u]
}$$

By  Fundamental Theorem \ref{2.5} and Lemma \ref{3.5}, the partial morphisms of geometries $\varphi_{x_i}\colon \mathbb{P}(V/\langle v_i\rangle)\dashrightarrow\mathbb{P}(V'/\langle v_i'\rangle)$ are projective, that is to say, they are the projectivization of some semilinear maps $$\begin{CD}
\Phi_i: V/\langle v_i\rangle @>>> V'/\langle v_i'\rangle\qquad (i=1,2)
\end{CD}$$ unique up to a scalar factor.

Since $\varphi_{x_i}(L/\!\!/x_i)=L'/\!\!/x_i'$ we have $\Phi_i(\langle v_1,v_2)/\langle v_i\rangle)\subseteq \langle v_1',v_2')/\langle v_i'\rangle$. Therefore, the maps $\Phi_i$ induce semilinear maps
$$\begin{CD}
\overline{\Phi}_i: V/\langle v_1,v_2\rangle @>>> V'/\langle v_1',v_2'\rangle\quad .
\end{CD}$$

Now, the projectivizations of $\overline{\Phi}_1$ and $\overline{\Phi}_2$ are just the partial morphism $\varphi_{L_0}$, so that  $\overline{\Phi}_1$ and $\overline{\Phi}_2$ differ in a scalar factor by \ref{1.14}. Multiplying $\Phi_1$ by this factor, we may assume that $\overline{\Phi}_1=\overline{\Phi}_2$. Therefore,  the respective morphisms $K\to K'$ associated to $\Phi_1$ and $\Phi_2$ coincide, since both also are associated to $\overline{\Phi}_1=\overline{\Phi}_2$.

Now, the equality $\overline{\Phi}_1=\overline{\Phi}_2:V/\langle v_1,v_2\rangle \longrightarrow V'/\langle v_1',v_2'\rangle$ let us define a semilinear map $\Phi\colon V\to V'$ via the diagram 
$$\begin{CD}
V @>{\Phi}>> V'\\
||| & & |||\\
V/\langle v_1\rangle\,\times_{(V/\langle v_1,v_2\rangle)}\,V/\langle v_2\rangle @>{\quad(\Phi_1,\Phi_2)\quad}>> V'/\langle v_1'\rangle\,\times_{(V'/\langle v_1',v_2'\rangle)}\,V'/\langle v_2'\rangle
\end{CD}$$
where the vertical isomorphisms are due to Lemma \ref{3.1}.

Let $\phi\colon\mathbb{P}\dashrightarrow\mathbb{P}'$ be the projectivization of $\Phi\colon V\to V'$. Let us see that this partial projective morphism is defined on $X$ and that it coincides with $\varphi$.

Given $x\in X$, $x'=\varphi(x)$, put $x=\mathbb{P}\langle v\rangle$, $x'=\mathbb{P}\langle v'\rangle$ for some vector subspaces $\langle v\rangle\subset V$, $\langle v'\rangle\subset V'$. We distinguish some cases according to the position of $x$.\smallskip

1. Case $x\notin\varphi^{-1}(x_1'\vee x_2')$, that is to say, $v'\notin\langle v_1',v_2'\rangle$.
The projective morphism $\varphi_{x_1}\colon\mathbb{P}/\!\!/x_1=X/\!\!/x_1\dashrightarrow X'/\!\!/x_1'=\mathbb{P}'/\!\!/x_1'$ fulfills $\varphi_{x_1}([x])=[\varphi(x)]=[x']$, so that the corresponding semilinear map  $\Phi_1\colon V/\langle v_1\rangle\to V'/\langle v_1'\rangle$ fulfills $\Phi_1([v])= [\lambda_1v']$ for some non null scalar $\lambda_1\in K'$. Analogously, $\Phi_2\colon V/\langle v_2\rangle\to V'/\langle v_2'\rangle$ fulfills $\Phi_2([v])= [\lambda_2v']$. Therefore, according to the definition of $\Phi$ we have $$\Phi(v)\equiv(\Phi_1,\Phi_2)([v],[v])=([\lambda_1v'],[\lambda_2v'])\equiv \lambda v'$$
where $\lambda_1=\lambda_2=:\lambda$ because $[\lambda_1v']=[\lambda_2v']\neq 0$ in $V'/\langle v_1',v_2'\rangle$ (recall that $v'\notin\langle v_1',v_2'\rangle$).
Hence $\,\phi(x)=x'=\varphi(x)$.\smallskip

2. Case $x\in\varphi^{-1}(x_1')$ (the case $x\in\varphi^{-1}(x_2')$ is similar). Condition $\varphi(x)=x_1'$ implies that the morphism $\varphi_{x_1}\colon\mathbb{P}/\!\!/x_1=X/\!\!/x_1\dashrightarrow\mathbb{P}\,'/\!\!/x_1'$ is not defined on $[x]$, hence $[v]\in\text{ker}\,\Phi_1$. On the other hand, the morphism $\varphi_{x_2}\colon\mathbb{P}/\!\!/x_2=X/\!\!/x_2\dashrightarrow\mathbb{P}\,'/\!\!/x_2'$ fulfills $\varphi_{x_2}([x])=[x_1']$, so that the associated semilinear map $\Phi_2\colon V/\langle v_2\rangle\to V'/\langle v_2'\rangle$ fulfills $\Phi_2([v])=[\lambda v_1']$ for some scalar $\lambda\neq 0$. Then
$$\Phi(v)\,\equiv\,(\Phi_1,\Phi_2)([v],[v])\,=\, (0,[\lambda v_1'])\,\equiv\, \lambda v_1'$$
and we conclude that $\phi(x)=x_1'\,=\,\varphi(x)$.\smallskip

3. Case $x\in\varphi^{-1}(x_1'\vee x_2')$, $x\not\in\varphi^{-1}(x_1')$, $x\notin\varphi^{-1}(x_2')$. That is to say, $x'\in x_1'\vee x_2'$, $x'\neq x_1'$, $x'\neq x_2'$, hence $v'=\alpha_1v_1'+\alpha_2v_2'\in\langle v_1',v_2'\rangle$ where $\alpha_1\neq 0$, $\alpha_2\neq 0$. The argument of case 1 shows that $\,\Phi_i([v])=[\lambda_iv']\,\in V'/\langle v_i'\rangle$, hence
$$\Phi(v)\equiv(\Phi_1,\Phi_2)([v],[v])=([\lambda_1v'],[\lambda_2v'])=([\lambda_1\alpha_2v_2'],[\lambda_2\alpha_1v_1'])$$
$$\equiv\, \lambda_2\alpha_1v_1'+\lambda_1\alpha_2v_2'$$
hence $\phi(x)\in x_1'\vee x_2'$.

Since $\text{dim}\,\varphi(X)\geq 3$, there exists a plane of $\varphi(X)$ intersecting the line $x_1'x_2'$ just at the point $x'$. This plane is generated by $x'$ and two other points $z_1',z_2'$. Let us consider a plane $z_1z_2x$ of $X$, where $z_1\in\varphi^{-1}(z_1')$, $z_2\in\varphi^{-1}(z_2')$. Since $\varphi\colon z_1z_2x\to z_1'z_2'x'$ is injective, the  plane $z_1z_2x$ intersects $\varphi^{-1}(x_1'x_2')$ just at the point $x$. Since $X$ has enough points, we may consider a different fourth point $z_3$ in the plane $z_1z_2x$ such that $z_1 z_2x=z_1 z_2 z_3$. Applying $\varphi$ to this equality one has $ z_1'\vee z_2'\vee x'=z_1'\vee z_2'\vee z_3'$. Therefore,
$$\phi(x)\,\in\,\phi(z_1 z_2x=z_1 z_2 z_3)\,\subseteq\,\phi(z_1)\vee\phi(z_2)\vee\phi(z_3)\,\overset{\text{case 1}}{=}\,z_1'\vee z_2'\vee z_3'\,=\,  z_1'\vee z_2'\vee x'$$

Since $ z_1'\vee z_2'\vee x'$ intersects $x_1'\vee x_2'$ just at $x'$, and $\phi(x)\in x_1'\vee x_2'$, we conclude that $\,\phi(x)=x'\,=\,\varphi(x)$.

\hfill$\square$\medskip

In \ref{5.3} we give a more general version of this theorem.

\begin{remark}\label{2.7} Under the hypothesis of the theorem, if $\varphi\colon X\to X'$ is injective then the extension $\phi\colon\mathbb{P}\dashrightarrow\mathbb{P}\,'$ is globally defined.
\end{remark}
{\it Proof.} If the exceptional subspace $E$ of $\phi\colon\mathbb{P}\dashrightarrow\mathbb{P}\,'$ is not empty, we consider the line $L$ of $\mathbb{P}$ joining $x_1$ and a point of $E$. Since $L\cap X$ is a line of $X$ (because it corresponds to a point of $\mathbb{P}/\!\!/x_1=X/\!\!/x_1$), it contains a point $\bar x_1\neq x_1$. Then $x_1\vee E=\bar x_1\vee E$, hence $\phi(x_1)=\phi(\bar x_1)$, that is to say, $\varphi(x_1)=\varphi(\bar x_1)$, against the injective character of $\varphi$.

\hfill$\square$

\begin{remark}\label{2.8} Under the hypothesis of the theorem, if $\varphi\colon X\to X'$ is an embedding then the extension $\phi\colon\mathbb{P}\to\mathbb{P}\,'$ is an embedding.
\end{remark}
{\it Proof.} First assume $\text{dim}\, X<\infty$. We have $\text{dim}\, X=\text{dim}\,\mathbb{P}$ since $\,\text{dim}\,X=\text{dim}(X/\!\!/x_1)+1=\text{dim}(\mathbb{P}/\!\!/x_1)+1=\text{dim}\,\mathbb{P}$. Then, from the inequalities $\text{dim}\,\mathbb{P}\geq\text{dim}\,\phi(\mathbb{P})\geq\text{dim}\,[\phi(X)=\varphi(X)]=\text{dim}\, X$, it follows that $\text{dim}\,\mathbb{P}=\text{dim}\,\phi(\mathbb{P})$, hence $\phi\colon\mathbb{P}\to\mathbb{P}\,'$ is an embedding by \ref{1.5a}b.

In the general case, we put $X=\bigcup X_i$, union of the directed family of subspaces $X_i$ of finite dimension, containing $x_1,x_2$, such that $\text{dim}\,\varphi(X_i)>2$. Let $\mathbb{P}_i$ be the projective closure of $X_i$ in $\mathbb{P}$. According to the former case, $\phi\colon\mathbb{P}_i\to\mathbb{P}\,'$ is an embedding, hence $\phi\colon\mathbb{P}=\bigcup\mathbb{P}_i\longrightarrow \mathbb{P}\,'$ is an embedding by \ref{1.6a}.

\hfill$\square$

\section{Locally projective geometries}

\noindent\textbf{Definition.} We say that a geometry $X$ fulfills the \textbf{local dimension formula} at a point $x\in X$ if for any pair of subspaces $S_1,S_2\subseteq X$ of finite dimension passing through $x$ we have
$$\text{dim}\,S_1+\text{dim}\, S_2\,=\, \text{dim}\,(S_1\vee S_2)+\text{dim}\, (S_1\cap S_2)\quad .$$

\begin{proposicion} Given a geometry $X$ and a point $x\in X$, the following statements are equivalent:\smallskip

a) $X$ fulfills the local dimension formula at  $x$.\smallskip

b) The quotient $X/\!\!/x$ is a projective space.
\end{proposicion}
{\it Proof.} ($a\Rightarrow b$). It is easy to check that for any pair of subspaces $S_1/\!\!/x$, $S_2/\!\!/x$ of $X/\!\!/x$ we have
$$(S_1/\!\!/x)\vee (S_2/\!\!/x)\,=\,(S_1\vee S_2)/\!\!/x\qquad ,\qquad (S_1/\!\!/x)\cap (S_2/\!\!/x)\,=\,(S_1\cap S_2)/\!\!/x\quad .$$
Moreover, $\text{dim}\,(S/\!\!/x)=\text{dim}\, S-1$.

Therefore,
$$\text{dim}\,\left((S_1/\!\!/x)\vee (S_2/\!\!/x)\right)+\text{dim}\,\left((S_1/\!\!/x)\cap (S_2/\!\!/x)\right)$$
$$=\,\text{dim}\,\left((S_1\vee S_2)/\!\!/x\right)+\text{dim}\,\left((S_1\cap S_2)/\!\!/x\right)$$
$$=\text{dim}\,(S_1\vee S_2)-1+\text{dim}\,(S_1\cap S_2)-1$$
(since $S_1,S_2$ meet at $x$, we may apply the local dimension formula)
$$=\,\text{dim}\,S_1+\text{dim}\,S_2-2\,=\,\text{dim}\,(S_1/\!\!/x)+\text{dim}\,(S_2/\!\!/x)$$

Hence in $X/\!\!/x$ the dimension formula holds. By \ref{1.10}, $X/\!\!/x$ is a projective space.\smallskip

($a\Leftarrow b$). Let $S_1,S_2$ be subspaces of $X$ of finite dimension such that $x\in S_1\cap S_2$. Since $S_1/\!\!/x$, $S_2/\!\!/x$ fulfill the dimension formula, the above argument shows that $S_1,S_2$ also do.

\hfill$\square$\medskip

For other equivalent properties see \cite{FF}, 4.4.1.\goodbreak

\begin{nada}\label{3.1a}\textbf{Definition.} A geometry $X$ is \textbf{locally projective} if $X/\!\!/x$ is a projective space for all $x\in X$. Equivalently, $X$ fulfills the local dimension formula at any point, that is to say, we have
$$\mathrm{dim}\,(S_1\vee S_2)+\mathrm{dim}\,(S_1\cap S_2)\,=\,\mathrm{dim}\,S_1+\mathrm{dim}\,S_2$$
for any pair of subspaces $S_1,S_2\subseteq X$ of finite dimension such that $S_1\cap S_2\neq\emptyset$.
\end{nada}

This notion only is relevant in dimension $\geq 3$, since any geometry of dimension  $\leq 2$ is locally projective.\goodbreak

\begin{nada}\label{3.3} Locally projective geometries (of dimension $\geq 3$) almost always have enough points and, therefore, the quotients $X/\!\!/x$ are irreducible projective spaces. Precisely: {\it If a locally projective geometry $X$ of dimension $\geq 3$ has not enough points then it is a reducible projective space} (see \cite{Wyler} Th. 6.2).
\end{nada}

\begin{proposicion}\label{3.3b} Any locally projective geometry $X$ is generated by lines and planes.
\end{proposicion}
{\it Proof.} Let $A\subseteq X$ be a subset containing the line joining any two different points of $A$ and the plane passing through any three non collinear points of $A$. Let us see that $A$ is a subspace. If $A$ is empty or a point, it is clear. Otherwise, fix $a\in A$, then $A=\bigcup L_i$ where $L_i$ runs over the lines of $X$ such that $a\in L_i\subseteq A$. Moreover, any two such lines $L_i\neq L_j$ fulfill that the plane $P_{ij}=L_i\vee L_j$ is contained in $A$.

Since $X/\!\!/a$ is generated by lines (because it is a projective space), we have that the set $[A]=\bigcup[L_i]$ is a subspace of $X/\!\!/a$, since any pair of points $[L_i],[L_j]$ is contained in the line $[P_{ij}=L_i\vee L_j]$. Therefore, $[A]=S/\!\!/a$ for some subspace $S$ of $X$. We have $\,S=\bigcup L_i=A$, hence $A$ is a subspace of $X$.

\hfill$\square$\medskip

From the above result and \ref{1.7} it follows the following

\begin{corolario} A map $\varphi\colon X\to X'$, between locally projective geometries, is a morphism if and only if for any $\,x_0,x_1,x_2,x_3\in X\,$ we have 
$$x_0\,\in\, x_1\vee x_2\vee x_3\qquad\Rightarrow\qquad\varphi(x_0)\,\in\,\varphi(x_1)\vee\varphi(x_2)\vee\varphi(x_3)\qquad .$$
\end{corolario}

Proposition \ref{3.3b} suggests that locally projective geometries may be defined by axioms of incidence of points, lines and planes. This is provided by the following alternative definition.

\begin{nada}\label{4.5} \textbf{Alternative definition.} A \textbf{locally projective geometry} is a set $X$ (whose elements are named \textbf{points}) with two families $\mathfrak{L}$ and $\mathfrak{P}$ of subsets (named \textbf{lines} and \textbf{planes} respectively) fulfilling the following four axioms:\smallskip

\textbf{LP1.} There is a unique line passing through any two different points and any line has at least two different points.\smallskip

\textbf{LP2.} There is a unique plane passing through any three non collinear points and any plane has at least three non collinear points.\smallskip

\textbf{LP3.} The line joining any two different points of a plane is contained in the plane.\smallskip

From the above axioms it follows that given a line $L$ and a point $x\notin L$ there is a unique plane passing through $L$ and $x$ (it is denoted by $L\vee x$).\smallskip

\textbf{LP4.} Given two different lines $L_1,L_2$ of a plane $P$ and a point $x\notin P$ then the intersection of the planes $L_1\vee x$, $L_2\vee x$ is a line.\goodbreak\medskip

In the context of this alternative definition, a \textbf{subspace} of $X$ is a subset $S$ containing the line passing through any two different points of $S$ and the plane passing through any three non collinear  points of $S$. \end{nada}

In \cite{Wyler} there is another axiomatization in terms of the incidence relation of points and planes.

In the projective geometry, the Veblen--Young axiom means that any two different lines of a plane intersect at a point, but the axiom is formulated avoiding the notion of plane. Analogously, axiom LP4 means that any two different planes of a $3$-space are disjoint or the intersection is a line (the intersection is not a point), but it is formulated avoiding the notion of $3$-space.\medskip

Locally projective geometries of dimension $3$ may be axiomatized replacing LP4 by the axiom \textbf{LP4'}: {\it If two different planes intersect, the intersection is a line}, and adding the axiom \textbf{LP5}: {\it There are four non coplanar points}. It is interesting to point that axioms LP1, LP2, LP3, LP4', LP5 exactly encode our intuition on the incidence relations of points, lines and planes in the three-dimensional physical space; in particular, they do not involve the parallelism, a notion which is not experimental since it has not a local character. %Solo cabría añadir a esos axiomas el hecho of que el espacio físico parece tener numerosos puntos, aunque quizás no infinitos.
\medskip

It is easy to check that any locally projective geometry in the sense of \ref{3.1a} fulfills the axioms of definition \ref{4.5}. The converse is not so easy and may be found in \cite{Wyler}.\bigskip

The following proposition let us to construct examples of locally projective geometries.

\begin{proposicion}\label{3.2} Let $\mathbb{P}$ be a projective space (or more generally a locally projective geometry) and let $X\subseteq \mathbb{P}$ be a subset with the following property:\smallskip

Any line $L$ of $\mathbb{P}$ fulfills $\,L\cap X=\emptyset\,$ or $\,|L\cap X|\geq 2$.\smallskip

\noindent Then the subgeometry $X$ is locally projective.
\end{proposicion}
{\it Proof.} It is enough to show that, for any point $x_0\in X$, the inclusion $X/\!\!/x_0\subseteq \mathbb{P}/\!\!/x_0$ is an equality. For all $[p]\in \mathbb{P}/\!\!/x_0$, the line $x_0\vee p$ of $\mathbb{P}$ intersects $X$ at least at two different points $x_0,x$, so that $x_0\vee p=x_0\vee x$, hence $[p]=[x]$ in $\mathbb{P}/\!\!/x_0$, that is to say, $[p]\in X/\!\!/x_0$.

\hfill$\square$

\begin{nada}\label{3.7d} \textbf{Examples}\smallskip

The following examples fulfill the hypothesis of the above proposition, so that they define locally projective geometries.\smallskip

\textbf{1.} {\it  Let $Z$ be a locally projective geometry whose lines have at least three points (for example, when $Z$ is an irreducible projective space). Then for any subspace $S\subset Z$ we have that $X=Z-S$ is locally projective.}

A particular case is an affine space, which is the complement of an hyperplane in an irreducible projective space.\smallskip

\textbf{2.} {\it Let $\mathbb{P}=\mathbb{P}(V)$ be a projective space over a division ring $K$ and let $\{ S_i\}$ be a family of subspaces of $\mathbb{P}$. If the  cardinal of the family $\{ S_i\}$ is $<|K|$, then $X=\mathbb{P}-\bigcup S_i$ is locally projective.}\smallskip

\textbf{3.} Let $\mathbb{P}$ be a projective space over a division ring $K$ of infinite order. In general, any open set $U$ of a  reasonable topology on $\mathbb{P}$ fulfills the hypothesis of Proposition \ref{3.2}, so that $U$ is locally projective. For example,\smallskip

{\it Let $\mathbb{P}_n$ be the projective space of finite dimension $n$ over $\mathbb{R}$ or $\mathbb{C}$, with the usual transcendent topology. Any open set $U$ of $\mathbb{P}_n$ is locally projective.}\smallskip

{\it Let $\mathbb{P}_n$ be the projective space of finite dimension $n$ over an infinite field $K$, with the Zariski topology. Any open set $U$ of $\mathbb{P}_n$ is locally projective.}\smallskip

\textbf{4.}  {\it Let $H_1\dots,H_{n+1}$ be hyperplanes, with empty intersection, of a projective space $\mathbb{P}_n$ of dimension $n$. Any subgeometry $X\subseteq\mathbb{P}_n$ containing these hyperplanes is locally projective.}\smallskip

\textbf{5.} {\it Let $H_1,H_2$ be two different hyperplanes of an irreducible projective space. Then $X=(H_1\cup H_2)-(H_1\cap H_2)$ is a locally projective geometry.}\smallskip

\textbf{6.} {\it Let $K\subset L$ be a strict inclusion of division rings and let $\,\mathbb{P}^n_K=\mathbb{P}(K^{n+1})\subset\mathbb{P}(L^{n+1})=\mathbb{P}^n_{L}\,$ be the corresponding projective spaces of dimension $n$. Then $X=\mathbb{P}^n_{L}-\mathbb{P}^n_{K}$ is a locally projective geometry.}
\end{nada}

\medskip

\noindent\textbf{Definition.} We say that a geometry \textbf{fulfills the bundle theorem} when its lines have the following property:\smallskip

{\it Let $L_1,L_2,L_3,L_4$ be four lines, so that no three of them are coplanar. If five of the six pairs $\{L_i,L_j\}_{i< j}$ are coplanar then so is the sixth pair.}\medskip

The following facts are elementary:\smallskip

-- Any projective space fulfills the bundle theorem. In this case the four lines $L_1,L_2,L_3,L_4$ are concurrent.\smallskip

-- If a geometry $X$ fulfills the bundle theorem then so does any subgeometry $Y\subseteq X$. Therefore, this is a necessary condition so that a geometry may be embedded in a projective space.\medskip

It is not difficult to prove that any locally projective geometry of dimension $\geq 4$ fulfills the bundle theorem (see \cite{Wyler}, Th 5.2). However there are locally projective geometries of dimension $3$ which do not fulfill the bundle theorem (see for example \cite{Kantor} \S 5 Ex. 5).\medskip

\begin{nada}\label{4.8}\textbf{Embedding Theorem}\smallskip

{\it a) Any locally projective geometry $X$ of dimension  $\geq 4$ admits an embedding $X\hookrightarrow\mathbb{P}$ into a projective space $\mathbb{P}$} (Wyler \cite{Wyler}; see also \cite{Kantor} and  \cite{FF}).\smallskip

{\it b) Any locally projective geometry $X$ of dimension $3$, fulfilling the bundle theorem, admits an embedding $X\hookrightarrow\mathbb{P}$ into a projective space} (Kahn \cite{Kahn}; see also  \cite{Kreuzer1}).\smallskip

{\it Moreover, these embeddings $X\hookrightarrow\mathbb{P}$ are minimal in the sense that $\,X/\!\!/x=\mathbb{P}/\!\!/x\,$ for all $x\in X$.}
\end{nada}

\begin{remark}\label{4.9} Let $X$ be a locally projective geometry fulfilling the hypotheses of the above theorem and let $X\hookrightarrow\mathbb{P}$ be the minimal projective embedding. Assume that $X$  is not a reducible projective space; then  $X$ has enough points (see \ref{3.3}) and, therefore, the projective spaces $\,X/\!\!/x=\mathbb{P}/\!\!/x$ are irreducible. Then the projective space $\mathbb{P}$ also is  irreducible and, since $\text{dim}\,\mathbb{P}=\text{dim}\,X\geq 3$, according to \ref{2.17} we have $\mathbb{P}=\mathbb{P}(V)$ for some vector space $V$.
\end{remark}

Now we may state the Fundamental Theorem for locally projective geometries:

\begin{teorema}\label{3.10} Let $X,X'$ be locally projective geometries fulfilling the hypotheses of \ref{4.8} and \ref{4.9}. Let $X\hookrightarrow\mathbb{P}(V)$, $X'\hookrightarrow\mathbb{P}(V')$ be the  respective minimal projective embeddings.

Any morphism of geometries $\varphi\colon X\to X'$, whose image is not contained in a plane, is the restriction of a unique partial projective morphism $\,\mathbb{P}(V)\dashrightarrow\mathbb{P}(V')$.
\end{teorema}
{\it Proof.} Just apply Theorem \ref{3.6} to the morphism $\varphi\colon X\to\mathbb{P}(V')$.

\hfill$\square$

\begin{remark} As we have indicated in Remarks \ref{2.7} and \ref{2.8}, if the morphism $X\to X'$ is injective then the extension $\mathbb{P}(V)\to\mathbb{P}(V')$ is globally defined, and if $X\to X'$ is an embedding then so is the extension $\mathbb{P}(V)\to\mathbb{P}(V')$. 
\end{remark}

\begin{remark} The condition on the image of $\varphi$ may be improved in some cases:\smallskip

\textbf{a.} {\it If $X$ is an open set of $\mathbb{P}(V)$, with respect to some linear topology, and $\varphi\colon X\to \mathbb{P}(V')\neq *$ is an embedding such that the image is not contained in a line, then $\varphi$ extends to an embedding $\mathbb{P}(V)\to\mathbb{P}(V')\,$} (Frank 
\cite{Frank}, Satz 1).\smallskip

\textbf{b.} Let $\mathbb{P}$ and $\mathbb{P}'$ be real projective spaces, with $\text{dim}\,\mathbb{P}<\infty$. Let $X\subseteq\mathbb{P}$ be an open set; it is a locally projective geometry generated by lines. {\it Any morphism $\varphi\colon X\to\mathbb{P}'$, such that the image is no contained in a line, is the restriction of a partial projective morphism $\mathbb{P}\dashrightarrow\mathbb{P}'\,$} (Lenz \cite{Lenz}, Hilfssatz 3). This result was generalized by Frank (\cite{Frank}, Satz 2) for projective spaces over division rings, endowed with a linear topology.\smallskip

\textbf{c.} The case when $X$ is an affine space is considered in Proposition \ref{4.1} below.  Morphisms of geometries $\varphi\colon\mathbb{A}\to\mathbb{A}'$ between real affine spaces, whose image is not contained in a line, are  affine morphisms. This result was generalized by W. Zick for affine spaces over arbitrary division rings (see \cite{Sancho} Th. 5.1).
\end{remark}

\section{Locally affino-projective geometries}

\noindent\textbf{Definition.} A subgeometry $X$ of an irreducible projective space $\mathbb{P}$ is said to be \textbf{affino-projective} if there is an hyperplane $H$ of $\mathbb{P}$ such that $X\cup H=\mathbb{P}$.

 Considering the affine space $\mathbb{A}=\mathbb{P}-H$, we have inclusions $\,\mathbb{A}\subseteq X\subseteq\mathbb{P}$.

\begin{proposicion}\label{4.1} Let $\mathbb{P}=\mathbb{P}(V)$ and $\,\mathbb{P}\,'=\mathbb{P}(V')\neq *$ be projective spaces over division rings $K$ and $K'$, respectively, with $|K|\geq 4$ or $|K|=3=\mathrm{char}\, K'$. Given an hyperplane $H\subset\mathbb{P}$, let us consider the affine space $\mathbb{A}=\mathbb{P}\!-\!H$.

Any morphism $\varphi\colon\mathbb{A}\to\mathbb{P}\,'$, whose image is not contained in a line, is the restriction of a unique partial projective morphism $\phi\colon\mathbb{P}\dashrightarrow\mathbb{P}'$.
\end{proposicion}

This proposition is proved in (\cite{Sancho}, Th. 4.10) and it is used to obtain a result of W. Zick determining the morphisms of geometries $\varphi\colon\mathbb{A}\to\mathbb{A}'$ between affine spaces. Let us extend the above proposition to affino-projective geometries.\goodbreak

\begin{proposicion}\label{5.2} With the hypotheses of \ref{4.1} on $\mathbb{P}$ and $\mathbb{P}\,'$, let $X\subseteq\mathbb{P}$ be an affino-projective subgeometry.
Any morphism $\varphi\colon X\to\mathbb{P}\,'$, whose image is not contained in a line, is the restriction of a unique partial projective morphism $\,\phi\colon\mathbb{P}\dashrightarrow\mathbb{P}'$.
\end{proposicion}
{\it Proof.} Let $H$ be an hyperplane of $\mathbb{P}$ such that $\,X\cup H=\mathbb{P}\,$ and put $\,\mathbb{A}:=\mathbb{P}-H\subseteq X$. Remark that $\varphi(\mathbb{A})$ is not contained in a line $L'$ of $\mathbb{P}\,'$: If $\varphi(\mathbb{A})\subseteq L'$ then $\varphi(X)\subseteq L'$ (because any point of $X$ is collinear with two points of $\mathbb{A}$), against the hypothesis on $\varphi$. By Proposition \ref{4.1}, the morphism $\varphi_{|\mathbb{A}}\colon\mathbb{A}\to\mathbb{P}'$ may be extended to a unique partial projective morphism $\phi\colon\mathbb{P}\dashrightarrow\mathbb{P}\,'$.

We have to show that $\phi_{|X}=\varphi$; we only have to check the equality on the points $x\in X\!-\!\mathbb{A}\subseteq H$. Let $p_1,p_2\in\mathbb{A}$ be points such that $\varphi(x),\varphi(p_1),\varphi(p_2)$ are not collinear (they exist), so that $\,\varphi(x)=(\varphi(p_1)\vee\varphi(x))\cap(\varphi(p_2)\vee\varphi(x))$. The lines $p_1\vee x$ and $p_2\vee x$ of $\mathbb{P}$ are different and they meet $\mathbb{A}$ at some points $q_1\neq p_1$ and $q_2\neq p_2$, respectively. Applying $\varphi$ to the equality $$x\,=\,(p_1\vee x)\cap(p_2\vee x)=(p_1\vee q_1)\cap(p_2\vee q_2)$$ we obtain
$$\varphi(x)\,=\,(\varphi(p_1)\vee\varphi(x))\cap(\varphi(p_2)\vee\varphi(x))\,=\,(\varphi(p_1)\vee\varphi(q_1))\cap(\varphi(p_2)\vee\varphi(q_2))$$
$$=\,(\phi(p_1)\vee\phi(q_1))\cap(\phi(p_2)\vee\phi(q_2))\,\owns\,\phi(x)$$

\hfill$\square$

%\begin{remark} The condition $|K|\geq 4$ or $|K|=3=\mathrm{char}\, K'$ in the above propositions may be eliminated imposing that $X$ has enough points and that the image of $\varphi$ is not contained in a  plane. Under these conditions, the result is a consequence of the main theorem \ref{3.6}. \end{remark}

\begin{lema}\label{4.6} Let $X$ be an affino-projective subgeometry of a projective space $\mathbb{P}$. Given a subspace $E\subseteq X$, we have that $X/\!\!/E$ is an affino-projective subgeometry of $\mathbb{P}/\!\!/\overline E$, where $\overline E$ is the closure of $E$ in $\mathbb{P}$.
\end{lema}
{\it Proof.} By Corollary \ref{2.3}, $X/\!\!/E$ is a subgeometry of $\mathbb{P}/\!\!/\overline E$. Let us see that it is affino-projective.

Let $H$ be an hyperplane of $\mathbb{P}$ such that $X\cup H=\mathbb{P}$. Let us show that if $E\nsubseteq H$ then $X/\!\!/E=\mathbb{P}/\!\!/\overline E$ (and we are done), and if $E\subseteq H$ then $(X/\!\!/E)\cup (H/\!\!/\overline E)=\mathbb{P}/\!\!/\overline E$. In the last case, if $\overline E=H$ then $H/\!\!/\overline E=\emptyset$, hence $X/\!\!/E=\mathbb{P}/\!\!/\overline E$ (and we are done) and if $\overline E\neq H$ then $X/\!\!/E$ is an  affino-projective subgeometry of $\mathbb{P}/\!\!/\overline E$ with associated hyperplane $H/\!\!/\overline E$.\smallskip

Let us prove that if $E\nsubseteq H$ (hence $\overline E\nsubseteq H$) then $X/\!\!/E=\mathbb{P}/\!\!/\overline E$. We have to show that any class $[p]\in\mathbb{P}/\!\!/\overline E$ is represented by an element $x\in X$. From the strict inclusions $\overline E\subset p\vee \overline E$, $(p\vee \overline E)\cap H\subset p\vee \overline E$, it follows that there exists $x\in p\vee \overline E$ such that $x\notin \overline E$, $x\notin(p\vee \overline E)\cap H$ (hence $x\notin H$ and then $x\in X$). By the exchange axiom, $p\vee \overline E=x\vee \overline E$, hence $[p]=[x]$.

Finally, from the equality $X\cup H=\mathbb{P}$ it directly follows that if $E\subseteq H$ (hence $\overline E\subseteq H$) then $(X/\!\!/E)\cup (H/\!\!/\overline E)=\mathbb{P}/\!\!/\overline E$.

\hfill$\square$\medskip

The following result is a variation of Theorem \ref{3.6}.

\begin{teorema}\label{5.3} Let $\mathbb{P}=\mathbb{P}(V)$ and $\,\mathbb{P}\,'=\mathbb{P}(V')$ be projective spaces of dimensions $\geq 3$ over division rings $K$ and $K'$, respectively, with $|K|\geq 4$ or $|K|=3=\mathrm{char}\, K'$. Let $X\subseteq\mathbb{P}$ be a subgeometry with enough points and let $\varphi\colon X\to\mathbb{P}'$, $x\mapsto x'$, a morphism such that:\smallskip

$1)$ The image of $\varphi$ is not contained in a plane of $\,\mathbb{P}\,'$.\smallskip

$2)$ There are points  $x_1,x_2\in X$ such that $x_1'\neq x_2'$ and $X/\!\!/x_i$ is an affino-projective subgeometry of $\,\mathbb{P}/\!\!/x_i$, where $i=1,2$.\smallskip

Then $\varphi\colon X\to\mathbb{P}\,'$ is the restriction of a unique partial projective morphism $\phi\colon\mathbb{P}\dashrightarrow\mathbb{P}\,'$.
\end{teorema}
{\it Proof.} Let $F_i:=\varphi^{-1}(x_i')$ and let $\overline F_i$ be the closure of $F_i$ in $\mathbb{P}$, where $i=1,2$.

According to Proposition \ref{2.1b}, the morphism $\varphi\colon X\to\mathbb{P}\,'$ induces partial morphisms
$$\begin{CD} \varphi_{x_i}: X/\!\!/x_i & \,--\!\to\, & X/\!\!/F_i @>{\varphi_{F_i}}>> \mathbb{P}\,'/\!\!/x_i'\qquad ,\qquad [x]\longmapsto[x']\quad .
\end{CD}$$
By Lemma \ref{3.5}, the image of $\varphi_{x_i}$ is not contained in a line of $\mathbb{P}'/\!\!/x_i'$. By Lemma \ref{4.6}, $X/\!\!/F_i$ is an affino-projective subgeometry of $\mathbb{P}/\!\!/\overline F_i$, so that applying  Proposition \ref{5.2} we obtain that $\varphi_{F_i}\colon X/\!\!/F_i\to\mathbb{P}'/\!\!/x_i'$ is the restriction of a partial projective morphism $\phi_{F_i}\colon\mathbb{P}/\!\!/\overline F_i\dashrightarrow\mathbb{P}'/\!\!/x_i'$. Therefore, the composition partial morphism
$$\begin{CD} \varphi_{x_i}: X/\!\!/x_i & \,--\!\to\, & X/\!\!/F_i @>{\varphi_{F_i}}>> \mathbb{P}\,'/\!\!/x_i'\end{CD}$$
 is the restriction of the partial projective morphism
 $$\begin{CD} \phi_{x_i}: \mathbb{P}/\!\!/x_i & \,\,--\!\to\,\, & \mathbb{P}/\!\!/\overline F_i & \,\,\overset{\phi_{F_i}}{--\!\to}\,\, & \mathbb{P}\,'/\!\!/x_i'\quad .\end{CD}$$
 
 The rest of the proof follows the lines of that of Theorem \ref{3.6}.

%Now the proof of the theorem is analogous to the proof given in Theorem \ref{3.6}.

\hfill$\square$\goodbreak\medskip

\noindent\textbf{Definition.} A subgeometry $X$ of an irreducible projective space $\mathbb{P}$ is said to be \textbf{locally affino-projective} when for any point $x\in X$ one has that $X/\!\!/x$ is an affino-projective subgeometry of $\mathbb{P}/\!\!/x$.\medskip

A non empty subset $X\subseteq\mathbb{P}$ is a locally affino-projective geometry if and only if for any point $x\in X$ there is an hyperplane $H_x$ of $\mathbb{P}$, passing through $x$, with the following property: For any line $L$ of $\mathbb{P}$ passing through $x$ we have $|L\cap X|\geq 2$ or $L\subseteq H_x$.

\begin{nada}\textbf{Examples}\smallskip

\textbf{1.} Let $\mathbb{P}_n$ be a projective space of finite dimension $n$ over an algebraically closed field $K$. Any algebraic hypersurface $X$ (defined by an homogeneous polynomial equation $q(x_0,\dots ,x_n)=0$ in some homogeneous coordinates of $\mathbb{P}_n$), non singular and of degree $m>1$,  is a locally affino-projective subgeometry.\smallskip

\textbf{2.} Let $\mathbb{P}_n$ be a projective space of finite dimension $n$ over a field $K$. Let $C$ be a quadric of $\mathbb{P}_n$ with non singular rational points. Then $X=C-\text{Sing}\,C$ is a locally affino-projective geometry.

As a particular case, in the real projective space $\mathbb{P}_3$, we obtain the following locally affino-projective geometries: $X=C$ where $C$ is an ellipsoid ({\it geometry of M\"obius}),  $X=C$ where $C$ is a non singular ruled quadric (\textit{geometry of Minkowski}), $X=C-\text{Sing}\, C$ where $C$ is a cone and $\text{Sing}\ C$ is the vertex (\textit{geometry of Laguerre}). These three geometries are named  \textit{planes of Benz}. A fourth geometry is $X=C-\text{Sing}\, C$ where $C$ is a pair of planes and $\text{Sing}\,C$ is the intersection line; moreover this geometry is locally projective.\smallskip

 \textbf{3.} The geometry of M\"obius may be generalized as follows. Let $\mathbb{P}$ be an irreducible projective space of dimension $\geq 3$. A non empty subset $X\subset\mathbb{P}$ is a \textbf{geometry of M\"obius} \cite{Maurer67} when for any point $x\in X$ the union of the lines $L$ such that $L\cap X=x$ is an hyperplane $H_x$ of $\mathbb{P}$ (that is to say, $X/\!\!/x=\mathbb{P}/\!\!/x-H_x/\!\!/x$).

Any geometry of M\"obius $X$  is \textbf{locally affine} (that is to say, $X/\!\!/x$ is an affine space for any $x\in X$). Conversely, by a theorem of M\"aurer \cite{Maurer}, any locally affine geometry $X$ of dimension $\geq 4$ is isomorphic to a geometry of M\"obius. This result extends to locally affine geometries of dimension $3$ fulfilling  the bundle theorem (Kahn \cite{Kahn1},\cite{Kahn0}).

A geometry of M\"obius $X\subset\mathbb{P}$ is said to be an \textbf{ovoid} when for any line $L$ of $\mathbb{P}$ we have $|L\cap X|\leq 2$.  \end{nada}

As a particular case of Theorem \ref{5.3} we get the following result.

\begin{nada} \textbf{Fundamental Theorem (locally affino-projective geometries)}\smallskip

{\it Let $\mathbb{P}=\mathbb{P}(V)$ and $\,\mathbb{P}\,'=\mathbb{P}(V')$ be projective spaces of dimensions $\geq 3$ over division rings $K$ and $K'$ respectively, with $|K|\geq 4$ or $|K|=3=\mathrm{char}\, K'$. Let $X\subseteq\mathbb{P}$ be a locally affino-projective subgeometry, with enough points.

Any morphism $\varphi\colon X\to\mathbb{P}\,'$, whose image is not contained in a plane,  is the restriction of a unique partial projective morphism $\phi\colon\mathbb{P}\dashrightarrow\mathbb{P}\,'$.}
\end{nada}

\begin{remark}\label{4.7}
Under the hypotheses of the theorem, assume that there is no point $p\in\mathbb{P}\!-\!X$ such that $(p\vee x)\cap X=x$ for all $x\in X$. Now, if the morphism $\varphi\colon X\to\mathbb{P}\,'$ is injective then the extension $\phi\colon\mathbb{P}\to\mathbb{P}\,'$ is globally defined. If moreover $\varphi\colon X\to\mathbb{P}\,'$ is an embedding then so is the extension $\phi\colon\mathbb{P}\to\mathbb{P}\,'$.\smallskip

The proofs are similar to the proofs of Remarks \ref{2.7} and \ref{2.8}. %No es difícil comprobar que la condición de inexistencia del punto $p$ se verifica en las geometries de los ejemplos  1 y 2 of \ref{4.5}. No sabemos si también la cumplen las del ejemplo 3.
\end{remark}

\begin{remark} If $X\subseteq\mathbb{P}$ is a M\"obius geometry over a division ring $K$, with $|K|>2$, then $X$ is a geometry generated by lines and planes. The proof is similar to the given in Proposition \ref{3.3b}, using that any affine space $\mathbb{A}$ over $K$ is generated by lines.
\end{remark}

\begin{remark} Let us consider an $n$-dimensional sphere $\mathbb{S}_n\subset\mathbb{R}^{n+1}=\mathbb{A}_{n+1}\subset\mathbb{P}_{n+1}$ as a M\"obius subgeometry (of dimension $n+1$) of the real projective space $\mathbb{P}_{n+1}$. Note that the planes of the subgeometry $\mathbb{S}_n$ are the circles. As a particular case of Theorem \ref{4.7}, we obtain the following result of \cite{Li-Wang}, Th. 1.2..\smallskip

{\it Let $\varphi\colon\mathbb{S}_n\to\mathbb{S}_n$ be an injective morphism (i.e., a circle-to-circle injection). If the image of $\varphi$ is not contained in a circle, then $\varphi$ is the restriction of a projectivity $\,\mathbb{P}_{n+1}\to\mathbb{P}_{n+1}\,$ (i.e., $\varphi$ is a M\"obius map).}
\end{remark}

%\section*{Declarations}

%\begin{itemize} \item Conflict of interest: We have no conflicts of interest to disclose. \item Funding: No funding. \item Authors' contributions: Not applicable since the article has only one author. \item Data Availability statement: Data sharing not applicable to this article as no datasets were generated or analyzed during the current study. \end{itemize}

\end{document}